\documentclass[11pt]{article}

\usepackage[latin1]{inputenc}

\usepackage{amsfonts,amsmath,amssymb,amsthm,graphicx,epsfig,float}
\usepackage[T1]{fontenc}

\usepackage[english,francais]{babel}

{
  
  \newtheorem*{theoreme*}{Th\'eor\`eme}

  \newtheorem{definition}{D\'efinition}
  
\newtheorem*{proposition*}{Proposition}
\theoremstyle{remark}
  \newtheorem*{remarque*}{Remarque}
}

\newcommand{\Cc}{\mathbb{C}}

\newcommand{\Pp}{\mathbb{P}}

\newcommand{\II}{\mathcal{I}}
\renewcommand {\epsilon}{\varepsilon}

\renewcommand {\leq}{\leqslant}
\renewcommand {\geq}{\geqslant}

\title{{\bf Sur la laminarité de certains courants}}
\author{Henry de Thélin}
\date{}

\begin{document}
\maketitle


\def\figurename{{Fig.}}%
\def\proofname{Preuve}
\def\contentsname{Sommaire}%







\begin{abstract}
Nous montrons qu'une suite de courbes analytiques lisses de la boule unité du plan complexe, dont le genre croît au plus comme l'aire, converge vers une lamination dans un sens faible.
\end{abstract}

\selectlanguage{english}
\begin{center}
{\bf{Laminarity of some currents}}
\end{center}

\begin{abstract}
We show that a sequence of smooth analytic curves of the unit ball of the complex plane, for which the genus is bounded by the area, converges to a lamination in a weak sense.
\end{abstract}

\selectlanguage{francais}

Mots-clefs: courants, laminarité, limites de courbes, théorie d'Ahlfors.\\
AMS: 32U40, 32H50.

\section*{{\bf Introduction}}

Dans cet article, on s'intéresse à des limites de courbes analytiques. La question est de savoir si celles-ci ont conservé un caractère analytique.\\
Soit $C_n$ une suite de courbes analytiques de la boule unité $B$ de $\Cc^2$. Si l'aire des $C_n$ reste bornée alors, quitte à extraire, $C_n$ converge vers une courbe analytique. C'est le théorème de Bishop (voir  \cite{B}).\\
Quand l'aire n'est plus uniformément majorée, on espère créer une lamination comme limite des $C_n$. Cependant, un excès de genre peut ôter tout caractère analytique à cette limite: dans l'exemple de Wermer (voir \cite{DS}), une suite de courbes dont le genre augmente plus vite que l'aire converge vers un compact ne contenant aucun disque holomorphe.\\
Notre lamination limite sera comprise dans un sens faible, celui de courant laminaire introduit par Eric Bedford, Michael Ljubich et John Smillie dans \cite{BLS}. Un $(1,1)$-courant positif est laminaire s'il s'écrit localement comme une intégrale de courants d'intégration sur une famille de disques, hors d'un ensemble négligeable (voir le paragraphe \ref{préliminaires} pour plus de détails).\\
L'objet de cet article est alors de démontrer:

\begin{theoreme*}
Soit $C_n$ une suite de courbes analytiques lisses de la boule unité $B$ de $\Cc^2$ (s'étendant un peu au-delà de $\overline{B}$).\\
On note $A_n$ l'aire de $C_n$, $G_n$ le genre de $C_n$ et on suppose que $T_n=\frac{[C_n]}{A_n}$ converge vers un $(1,1)$-courant positif fermé $T$ de $B$ (toujours possible quitte à extraire).\\
Alors, si $G_n=O(A_n)$, $T$ est laminaire.
\end{theoreme*}

Cet énoncé est une version locale de résultats précédents, de nature globale dans $\Pp^2(\Cc)$: dans \cite{BS}, Bedford et Smillie montrent que, pour une suite de courbes rationnelles $C_n$ de $\Pp^2(\Cc)$, les limites de $\frac{[C_n]}{A_n}$ sont laminaires. Dans \cite{D}, leur résultat a été étendu par Romain Dujardin aux courbes algébriques de  $\Pp^2(\Cc)$ de genre en $O(A_n)$ à singularités raisonnables.\\
L'ingrédient principal, dans leur situation, est la formule de Riemann-Hurwitz. Dans notre cas, on ne peut pas l'utiliser faute de revêtements au-dessus des directions complexes. Elle est remplacée ici (comme dans \cite{BLS} et \cite{C}) par la théorie d'Ahlfors qui conduit à une inégalité de Riemann-Hurwitz approchée.\\
Signalons enfin l'intérêt du théorème pour l'étude de courants limites de  $\frac{[C_n]}{A_n}$ où $C_n$ est une courbe algébrique lisse de $\Pp^2(\Cc)$ (par exemple $C_n=f^{-n}(L)$ où $f$ est un endomorphisme holomorphe de $\Pp^2(\Cc)$ et $L$ une droite projective générique). En effet, si on sait trouver un ouvert où le genre se concentre peu (en $O(A_n)$), on en déduit la laminarité des courants limites dans celui-ci. On peut par exemple appliquer ce raisonnement pour montrer que le courant de Green $T$ associé à un endomorphisme holomorphe critiquement fini de $\Pp^2(\Cc)$ est laminaire en dehors du support de $T \wedge T$. Nous espérons, par la suite, montrer un tel résultat pour tous les endomorphismes holomorphes de $\Pp^2(\Cc)$.\\
\\
Voici le plan de ce texte: dans le premier paragraphe, on donnera le schéma de la preuve tandis que dans le second on  détaillera la démonstration du théorème.\\
\\
{\bf Remerciement:} Je tiens à remercier mon directeur de thèse Julien Duval pour son aide précieuse dans l'élaboration de cet article.

\section{{\bf Préliminaires }\label{préliminaires}}

\subsection{{\bf Courants laminaires }}

Notre référence pour la laminarité est l'article de Bedford, Ljubich et Smillie (voir \cite{BLS}).

\begin{definition}{\bf }

Un courant $T$ est {\bf laminé} dans un bidisque, si dans celui-ci il s'écrit $\int_{\tau}[\Delta_t] d\lambda(t)$; ici les $\Delta_t$ sont des graphes, deux à deux disjoints, de fonctions holomorphes au-dessus d'une des directions du bidisque, et $\lambda$ une mesure positive portée par une transversale aux graphes.

\end{definition}

Les courants laminés ont de bonnes propriétés de compacité, par le théorème de Montel:

\begin{proposition*}
Soit $T_n$ une suite de courants laminés dans un bidisque. Si la masse des mesures transverses $\lambda_n$ reste bornée, alors, quitte à extraire, $T_n$ converge vers un courant $T$ qui est laminé dans le bidisque.

\end{proposition*}

On pourrait imaginer une définition des courants laminaires comme étant des courants localement laminés. Dans notre situation cette définition serait trop forte: cela nous conduit à une notion  plus souple (dite laminarité faible dans \cite{BLS}):

\begin{definition}

Un courant $T$ est {\bf laminaire}, s'il s'écrit comme une somme $\sum{T_j}$ avec $T_j$ laminé dans $U_j$, et une compatibilité entre les $T_j$: quand deux graphes se rencontrent, leur intersection est un disque.

\end{definition}

La laminarité d'un courant est liée à la présence de graphes de fonctions holomorphes au-dessus d'une direction complexe dans son support.\\
Pour montrer qu'une limite de $T_n=\frac{[C_n]}{A_n}$ est laminaire, la méthode sera donc de construire des graphes dans les courbes $C_n$ au-dessus d'une direction complexe, puis de passer à la limite.

\subsection {{\bf Schéma de la preuve}\label{Schéma}}

\begin{figure}[c,h]
\begin{center}

\setlength{\unitlength}{1cm}
\begin{picture}(5,3)
\includegraphics{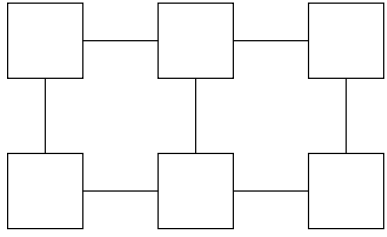}
\end{picture}

\end{center}
\caption{Pavage en croix{\label{Figure1}}. Les carrés font partie de la famille $Q$. Les croix pavent $C-Q$.}
\end{figure}

Dans ce paragraphe, on va montrer comment construire de bons disques sur les courbes $C_n$ dans un cas modèle.\\
On commence par se ramener à des courbes $C_n$ à géométrie bornée par $A_n$:\\ $L_n+G_n+B_n=O(A_n)$, où:
$$L_n=\mbox{longueur du bord de }C_n,$$
$$G_n=\mbox{genre de }C_n,$$
$$B_n=\mbox{nombre de composantes de bord de }C_n.$$
Puis, on se fixe une direction $D$ (qui vérifie $\pi_{*}T \neq 0$, où $\pi$ est la projection orthogonale associée à $D$), et on quadrille le carré $C \subset D$, centré en $0$, de côté $2$, en $4k^2$ carrés égaux. Un tel quadrillage peut être décomposé en quatre familles de $k^2$ carrés deux à deux disjoints.\\
On va partir de la famille $Q$ la moins recouverte. Elle vérifie en particulier:
 $$S_n(Q)=\frac{1}{\mbox{aire de } Q} \int_{C_n \cap \pi^{-1}(Q)}\pi^{*}\omega \leq S_n=\frac{1}{\mbox{aire de } C} \int_{C_n}\pi^{*} \omega$$
où $\omega$ est la forme kählérienne standard de $\Cc$.\\
Le but est de démontrer que les composantes connexes de $\pi^{-1}(Q)$ sont majoritairement des {\bf{îles}} (disques dont le bord se projette sur le bord des carrés de $Q$).\\
Le nombre d'îles au-dessus de $Q$ est lié à la caractéristique d'Euler de $C_n-\pi^{-1}(Q)$. En effet, si $\II$ désigne l'ensemble des îles de $\pi^{-1}(Q)$, on a $\chi(C_n - \pi^{-1}(Q)) \geq \chi(C_n) -\#\II$ (enlever une île fait chuter la caractéristique d'Euler de $1$).
En utilisant alors l'hypothèse sur le genre et le nombre de composantes de bord, on obtient une minoration du nombre d'îles par  $- \chi(C_n - \pi^{-1}(Q)) -O(A_n)$.\\
On va majorer  $\chi(C_n - \pi^{-1}(Q))$ pour obtenir une bonne minoration du cardinal de $\II$. Pour cela, on pave $C-Q$ en croix (voir figure \ref{Figure1}).

A partir de là, on construit un graphe où chaque sommet représente une composante connexe au-dessus d'une croix, et où on met une arête entre deux sommets si les composantes en question sont adjacentes. Dans toute la suite, on identifiera sommets et composantes connexes correspondantes.\\
On obtient alors:
\begin{equation*}
\begin{split}
\chi(C_n - \pi^{-1}(Q))&=\displaystyle\sum_{ \mbox{ sommets}}{\chi(\Sigma)}-\mbox{nombre d'arêtes}\\
&\leq s-a
\end{split}
\end{equation*}

où $s$ est le nombre de sommets et $a$ le nombre d'arêtes.\\
La combinaison de cette relation avec la précédente, nous conduit à une minoration du nombre d'îles par $a-s-O(A_n)$. Il nous reste donc à majorer le nombre de sommets et à minorer le nombre d'arêtes.\\
On va ajouter une hypothèse supplémentaire dans ce paragraphe: nous supposerons que les composantes connexes au-dessus des croix sont des graphes au-dessus de celles-ci.\\
Le nombre de sommets est alors égal à $k^2S_n(C-Q)$ où 
$$S_n(C-Q)=\frac{1}{\mbox{aire de }C-Q} \int_{C_n \cap \pi^{-1}(C-Q)}\pi^{*}\omega.$$
Pour minorer le nombre d'arêtes, on utilise la remarque suivante: si $\Sigma$ est un sommet et $\nu(\Sigma)$ sa valence (i.e. le nombre d'arêtes qui partent du sommet), on a:
$$L_n \geq \displaystyle\sum_{ \mbox{ sommets}}{(4-\nu(\Sigma))} \frac{1}{k},$$
c'est-à-dire:
$$\displaystyle\sum_{ \mbox{ sommets}}{\nu(\Sigma)} \geq 4k^2S_n(C-Q)-kO(A_n),$$
ce qui implique:
$$a \geq 2k^2S_n(C-Q)-kO(A_n).$$
Le nombre d'îles au-dessus de $Q$ est alors minoré par:
$$k^2S_n-kO(A_n)=k^2S_n(1-\epsilon),$$
car la direction $D$ est bien choisie ($\pi_{*}T \neq 0$).\\
Les composantes au-dessus de $Q$ sont donc bien majoritairement des îles.\\
Dans la réalité on n'aura pas vraiment des graphes au-dessus des croix, c'est pour cela que l'on utilisera un peu de théorie d'Ahlfors pour minorer le nombre d'arêtes.

\subsection{{\bf Un peu de théorie d'Ahlfors }\label{Ahlfors}}

Soit $f:\Sigma \longmapsto \Sigma_0$ une application holomorphe entre surfaces de Riemann, et $\gamma$ une courbe régulière de $\Sigma_0$. Un des premiers résultats de la théorie d'Ahlfors permet de comparer le nombre moyen de feuillets de $\Sigma$ au-dessus de $\Sigma_0$ avec le nombre moyen de feuillets au-dessus de $\gamma$. Plus précisément:\\
$$S=\frac{\mbox{aire de }f(\Sigma) \mbox{ comptée avec multiplicité }}{\mbox{aire de }\Sigma_0},$$
$$L=\mbox{la longueur du bord relatif de }\Sigma \mbox{ (i.e. la longueur de }f(\partial \Sigma)- \partial \Sigma_0),$$
$$S(\gamma)=\frac{\mbox{longueur de }f(f^{-1}(\gamma)) \mbox{ comptée avec multiplicité }}{\mbox{longueur de }\gamma}.$$
Alors, on a (voir par exemple \cite{N}):

\begin{theoreme*}
Il existe une constante $h$ indépendante de $f$ telle que:
$$|S-S(\gamma)| \leq hL.$$
\end{theoreme*}

\section {{\bf Démonstration du théorème}} 

\subsection{{\bf Simplification géométrique des courbes $C_n$}}

Voici comment on se ramène à des courbes $C_n$ à géométrie bornée par $A_n$, quitte à rétrécir un peu la boule de départ:\\
Considérons les trois boules concentriques $\rho B$, $\rho' B$ et $B$ (où $\rho'$ est strictement compris entre $\rho$ et $1$). Notons $\widetilde{C_n}$ la courbe obtenue en collant les composantes connexes de $C_n \cap(\rho' B - \rho B)$ qui touchent $\rho \partial B$, à $C_n \cap (\rho B)$. Par le principe du maximum, $\widetilde{C_n}$ a son bord inclus dans  $\rho' \partial B$. Nous allons voir que $\widetilde{C_n}$ a un nombre de composantes de bord en $O(A_n)$. En effet, si on coupe $C_n$ suivant ces $N$ composantes de bord, grâce au contrôle du genre on obtient au moins $N-O(A_n)$ composantes connexes $\Sigma$ dans l'une ou l'autre des calottes sphériques $B - \rho' B$, $\rho' B - \rho B$. Par construction le bord de $\Sigma$ doit rencontrer les deux sphères bordant la calotte dans laquelle elle se trouve. Un argument d'aire s'appuyant sur le théorème de Lelong (voir \cite{L}) montre alors que $N-O(A_n)$ est un $O(A_n)$, donc $N$ aussi. De plus, quitte à bouger un peu $\rho' \partial B$ et à extraire une sous-suite, on a  $\widetilde{L_n}=O(A_n)$ via la formule de coaire (voir \cite{M}).\\
On s'est donc ramené à une courbe $\widetilde{C_n}$ qui vérifie  $\widetilde{L_n}+\widetilde{G_n}+\widetilde{B_n}=O(A_n)$ et qui coïncide avec $C_n$ sur $\rho B$.\\
Dans la suite, tous les tildes seront oubliés.

\subsection{{\bf Mise en place de la preuve}}  

Reprenons la démonstration du paragraphe \ref{Schéma}:\\
On part d'une direction $D$ (qui vérifie $\pi_{*}T \neq 0$, où $\pi$ est la projection orthogonale associée à $D$). On quadrille le carré $C \subset D$, centré en $0$, de côté $2$ en $4k^2$ carrés égaux (associé à un tel quadrillage, il y a quatre familles de $k^2$ carrés deux à deux disjoints). On part de la famille $Q$ la moins recouverte (en particulier le nombre moyen de feuillets $S_n(Q)$ au-dessus de $Q$ est inférieur au nombre moyen de feuillets $S_n$ au-dessus du quadrillage initial).
Puis on pave $C - Q$ en croix comme dans le paragraphe précédent (voir figure \ref{Figure1}).




Dans le paragraphe \ref{Schéma}, on s'était placé dans un modèle où les composantes au-dessus des croix étaient des graphes. Ici ce n'est plus vrai, mais on va voir que l'on peut se ramener au cas où la majorité de ces composantes se projettent sur presque toute la croix correspondante.\\
Si $\Sigma$ est une composante au-dessus d'une croix, on a plusieurs possibilités. En effet soit $\epsilon_k$ une suite qui tend lentement vers $0$ (dans ce texte toute suite tendant vers $0$ sera notée $\epsilon_k$).\\
Si on note $l(\Sigma)$ la longueur du bord relatif de $\Sigma$ et $a(\Sigma)$ l'aire de $\pi(\Sigma)$ comptée avec multiplicité, on peut avoir:
$$l(\Sigma) \geq \frac{\epsilon_k}{k}$$ (on dira que $\Sigma$ a un bord long), ou:
$$l(\Sigma) \leq \frac{\epsilon_k}{k}$$ (on parlera de bord court).\\
En utilisant l'inégalité isopérimétrique on remarque que, dans le deuxième cas, on a soit:
$$a(\Sigma) \geq (1-\epsilon_k)\mbox{aire de la croix},$$
soit:
$$a(\Sigma) \leq 4 \mbox{ }l(\Sigma)^2.$$

 Dans notre contexte, ce sont les composantes du dernier type qui sont les plus éloignées de la situation du paragraphe \ref{Schéma}: on veut donc les enlever.\\
En ôtant ces composantes $(c_1,...,c_m)$, on modifie la courbe $C_n$ d'une aire (pour $\pi^* \omega$) au plus égale à $\displaystyle\sum_{i=1}^{m}{4l_i^2} \leq \frac{4}{k} \displaystyle\sum_{i=1}^{m}{l_i}$. D'après l'hypothèse sur $L_n$, on voit que cette modification est inférieure à $\frac{1}{k}O(A_n)$ (donc négligeable). D'autre part, la courbe obtenue en enlevant ces composantes vérifie toujours $L_n+G_n+B_n=O(A_n)$.\\
Dans la suite, on travaillera avec $C_n$ privée de ces composantes, que l'on notera toujours $C_n$.\\
\\
Notons $\II$ l'ensemble des îles au-dessus de $Q$. Dans le paragraphe \ref{Schéma} on a vu, via la construction du graphe où chaque sommet représente une composante connexe au-dessus d'une croix, et où l'on met autant d'arêtes entre deux sommets qu'il y a d'arcs en commun dans le bord des composantes correspondantes, que le cardinal de $\II$ est minoré par $a-s-O(A_n)$ (où $s$ est le nombre de sommets et $a$ le nombre d'arêtes) .\\
Il nous reste à majorer le nombre de sommets et à minorer le nombre d'arêtes.

\subsection{{\bf Majoration du nombre de sommets et minoration du nombre d'arêtes}}

Grâce à notre simplification effectuée au paragraphe précédent, on sait qu'un sommet a soit un bord long, soit un bord court et une aire pour $\pi^{*} \omega$ supérieure à  $(1-\epsilon_k)\mbox{aire de la croix}$.\\
D'après l'hypothèse sur la longueur du bord, il y a au plus $\frac{k}{\epsilon_k}O(A_n)$ sommets avec un bord long.\\
Pour les sommets avec une aire supérieure à $(1-\epsilon_k)\mbox{aire de la croix}$, on voit facilement qu'il y en a au plus $S_n(C - Q) \frac{k^2}{(1-\epsilon_k)}$ où $S_n(C - Q)$ est le nombre moyen de feuillets au-dessus de $C-Q$. On obtient alors une majoration de $s$ par: 
$$S_n(C - Q)k^2(1+\epsilon_k)+ \frac{k}{\epsilon_k}O(A_n).$$
Passons maintenant à la minoration du nombre d'arêtes, par la théorie d'Ahlfors (voir le paragraphe \ref{Ahlfors}).\\
On se fixe un sommet $\Sigma$ au-dessus d'une croix. Les côtés de la croix qui ne sont pas dans le bord de $Q$ seront notés $\alpha_1$, $\alpha_2$ , $\alpha_3$ et $\alpha_4$.\\
En appliquant la théorie d'Ahlfors à $f=k\pi:\Sigma \longmapsto \Sigma_0$, où $\Sigma_0$ est une croix dont la taille ne dépend plus de $k$, on obtient:
$$|S(\alpha_i)-S(\Sigma)| \leq hk \mbox{ }l(\Sigma),$$
où $h$ est indépendante de $k$, $S(\alpha_i)$ est le nombre moyen de feuillets de $\pi_{|\Sigma}$ au-dessus de $\alpha_i$ et $S(\Sigma)$ le nombre moyen de feuillets de $\pi_{|\Sigma}$ au-dessus de la croix. D'où,
$$\displaystyle\sum_{i=1}^4{S(\alpha_i)} \geq 4S(\Sigma)-4hk\mbox{ }l(\Sigma).$$
Comme on peut supposer, quitte à bouger le quadrillage initial, que les $\alpha_i$ ne contiennent pas de valeurs critiques de $\pi_{|\Sigma}$, le nombre de composantes des préimages $\pi^{-1}(\alpha_i)$ dans $\Sigma$ est supérieur à $\displaystyle\sum_{i=1}^4{S(\alpha_i)}$. On a donc: $\nu(\Sigma) \geq 4S(\Sigma)-4hk \mbox{ }l(\Sigma)$, où $\nu(\Sigma)$ est la valence du sommet $\Sigma$ dans le graphe (voir paragraphe \ref{Schéma}).\\
En considérant maintenant tous les sommets au-dessus de la croix, puis toutes les croix, on obtient:
$$2a=\displaystyle\sum_{ \mbox{ sommets}}{\nu(\Sigma)} \geq 4k^2S_n(C-Q)-kO(A_n),$$
d'où une minoration de $a$ par $2k^2S_n(C-Q)-kO(A_n)$, car la direction $D$ est bien choisie ($\pi_*T \neq 0$).\\
On a trouvé une minoration du nombre d'îles au-dessus de $Q$ en $S_nk^2(1-\epsilon_k)$ \\
(car $S_n(C - Q) \geq S_n$).\\\\
L'estimée précédente implique que $Q$ était quand même bien recouverte. En particulier, si $Q'$ est une famille de carrés, on a $S_n(C - Q') \geq (1-\epsilon_k)S_n$, d'où par ce que l'on a fait, il y a au moins $(1-\epsilon_k)k^2S_n$ îles au-dessus de $Q'$, soit $4(1-\epsilon_k)k^2S_n$ îles au-dessus du quadrillage initial.\\
De plus un argument d'aire montre que très peu d'entre elles sont ramifiées. On a donc au moins $4(1-\epsilon_k)k^2S_n$ {\bf{bonnes îles}} (graphes au-dessus des carrés du quadrillage) dans les courbes $C_n$.\\
L'objet du paragraphe suivant est d'en déduire la laminarité de $T$.

\subsection{{\bf Laminarité de $T$ dans $B$}}  
 
Soit $T_{k,n}$ le courant défini par $T_{k,n}=\frac{1}{A_n} \displaystyle\sum_{ \mbox{ bonnes îles}}[\Gamma]$     ($T_{k,n}$ est laminé au-dessus de chaque carré du quadrillage).\\
Rappelons que $T_n=\frac{[C_n]}{A_n}$.\\
La minoration du paragraphe précédent nous conduit à:
$$\int{T_{k,n} \wedge \pi^*\omega} \geq (1-\epsilon_k)\int{T_{n} \wedge \pi^*\omega},$$
d'où,
$$\int{(T_{n}-T_{k,n}) \wedge \pi^*\omega} \leq \epsilon_k.$$
En utilisant maintenant la proposition de compacité du paragraphe \ref{préliminaires}, $T_{k,n}$ converge vers un courant $T_k$ laminé au-dessus des carrés du quadrillage (quitte à extraire une sous-suite), et on a toujours l'estimée:
$$\int{(T-T_{k}) \wedge \pi^*\omega} \leq \epsilon_k,$$
avec $T-T_k \geq 0$ par construction.\\
Si on raffine de plus en plus le quadrillage (i.e. si $k$ augmente), $T_k$ croît vers un courant $T_{\infty}$ qui est laminaire (passer de $T_k$ à $T_{k'}$ avec $k' > k$ revient à rajouter des courants laminés qui sont compatibles entre eux). De plus, $T_{\infty} \leq T$ et $\int{(T-T_{\infty}) \wedge \pi^*\omega} \leq 0$.\\
Alors, si la direction $D$ est choisie de sorte que $T$ ne charge pas globalement les fibres de $\pi$, cela permet de conclure que $T=T_{\infty}$, c'est-à-dire que $T$ est laminaire dans $B$.

\bigskip

Henry de Thélin\\
Université Paul Sabatier, Laboratoire Emile Picard/CNRS UMR 5580, 118, route de Narbonne, 31062 Toulouse cedex 4, France.\\
E-mail: dethelin@picard.ups-tlse.fr

\end{document}